\documentclass[reqno,11pt]{amsart}

\usepackage{amsmath,amsfonts,amssymb,amsthm,epsfig,enumerate}
\usepackage{xcolor}
\usepackage[english]{babel}
\usepackage{hyperref} 
 \allowdisplaybreaks
\numberwithin{equation}{section}

\font\script=rsfs10 at 11pt
\def\eps{\varepsilon}
\def\H{{\mbox{\script H}\,\,}}

\def\L{{\mbox{\script L}\,\,}}
\def\E{\mathcal E}
\def\F{\mathcal F}
\def\G{\mathcal G}

\def\J{\mathcal J}

\def\R{\mathbb R}
\def\freccia#1{\xrightarrow[\ #1]{}}
\def\C{\mathcal C}
\def\S{\mathbb S}
\def\N{\mathbb N}
\def\p1{{P}}

\def\tv#1{\big| |#1|\big|}

\def\bal{\begin{aligned}}
\def\eal{\end{aligned}}
\def\proofof#1{\begin{proof}[Proof of #1]}
\def\pr#1#2{\langle #1,\,#2\rangle}
\def\Chi#1{\hbox{{\large $\chi$}{\Large $_{_{#1}}$}}}
\def\case#1#2{\par\noindent{\underline{\it Case~#1.}}\emph{ #2}\\}
\def\step#1#2{\par\noindent{\underline{\it Step~#1.}}\emph{ #2}\\}
\def\comp{\subset\subset}

\newcounter{mt}
\def\maintheorem#1#2#3{\par \medskip \noindent {\bf Theorem~\mref{#1}}~(#2).~{\it #3}\par}
\def\mref#1{\Alph{#1}}
\def\maintheoremdeclaration#1{\stepcounter{mt}\newcounter{#1}\setcounter{#1}{\arabic{mt}}}

\maintheoremdeclaration{main}

\newtheorem{theorem}{Theorem}[section]
\newtheorem{lemma}[theorem]{Lemma}
\newtheorem{prop}[theorem]{Proposition}

\newtheorem{defin}[theorem]{Definition}

\begin{document}

\title[Disconnected minimal clusters in Finsler manifolds]{Non-connected minimal clusters of small volume in a compact Finsler manifold}

\author{D. Carazzato}
\address{Davide Carazzato\\Faculty of Mathematics, University of Vienna (Austria)}
\email{davide.carazzato@unive.ac.at}

\author{A. Pratelli}
\address{Aldo Pratelli\\ Department of Mathematics, University of Pisa (Italy)}
\email{aldo.pratelli@unipi.it}

\begin{abstract}
We prove that isoperimetric clusters in general Finsler manifolds can be disconnected even when the volume constraint is small. This is in contrast with the Riemannian case and the ``fixed-norm'' case, where connectedness for arbitrarily small volume has recently been proved.
\end{abstract}

\maketitle

\section{Introduction}

In this paper, we are interested in the isoperimetric problem for clusters in a compact manifold $M$, and in particular we focus on the case of ``small'' minimal clusters, that is, minimal clusters with a very small total volume. Intuitively, since a manifold locally ``looks like'' $\R^N$, in some sense, one can think that such clusters should be more or less similar to corresponding clusters in $\R^N$. In particular, since minimal clusters in $\R^N$ are well-known to be connected, one might wonder whether sufficiently small minimal clusters in a compact manifold must be connected. For sets, which correspond to clusters with just one chamber, the behaviour of small isoperimetric regions has been widely investigated in several works: see for example~\cite{MJ2000,N2009,FM2011,DeRN2025} and the works related to the more sigular setting~\cite{ABFP2022,CM2022,AG2024,CM2024,APPS2025}, together with some surveys~\cite{P2023,A2025} and the references therein.\par

In the case of clusters, this kind of study is relatively recent, and a positive answer has been given in~\cite{NP2025}. More precisely, it has been shown that sufficiently small minimal clusters are connected if $M$ is a compact Riemannian manifold (the precise statement is given in Theorem~\ref{thmknown} below).

It is then natural to consider the more general case of Finsler manifolds. The formal definition of a Finsler manifold with some relevant bibliography will be given in Subsection~\ref{secFinsler}, but it is easy to get a rough idea. Namely, while the geometry around a point in a Riemannian manifold is basically that of $\R^N$ with the Euclidean norm (up to a linear change of coordinates), the geometry around a point in a Finsler manifold is basically that of $\R^N$ with some given convex norm, which is not necessarily the Euclidean one and which might depend on the point. As already pointed out in the paper~\cite{NP2025} itself, the proof in the Riemannian setting seems to rely mainly on the fact that the geometry around any point is somehow always the same, and not that this is precisely the Euclidean one. As a consequence, it was conjectured in~\cite{NP2025} that the same statement should hold true also for ``fixed-norm manifolds'', that is, manifolds for which the ``infinitesimal norm'' at any point is the same, up to a linear change of coordinates with determinant $\pm1$, see again Subsection~\ref{secFinsler} for the formal definition. The conjecture has been recently proved in~\cite{FPP}. Hence, these are the known results.

\begin{theorem}[\cite{NP2025,FPP}]\label{thmknown} Let $M$ be a compact (smooth or ${\rm C}^2$) Finsler manifold of fixed-norm type, and let $k\in\N$. Then, there exists some $\bar\eps>0$ such that any minimal $k$-cluster $\E\subseteq M$ with total volume $\tv{\E}<\bar\eps$ is connected.
\end{theorem}

As a byproduct of the above results, in~\cite{NP2025,FPP} it has been proved that the \emph{isoperimetric profile} of the manifold is asymptotically close to that of $\R^N$ endowed with the infinitesimal norm. Let us be more precise. The isoperimetric profile of the manifold $M$ for $k$-clusters is defined as
\[
\J_{M,k}(v) = \inf \Big\{ P_M(\E):\, \E\subseteq M,\, |\E|=v \Big\}\,,
\]
that is, for every admissible volume $v\in (\R^+)^k$, $\J_{M,k}(v)$ is the minimal perimeter of the $k$-clusters with volume $v$. Similarly, if $\|\cdot\|$ is any norm in $\R^N$, the isoperimetric profile $\J_{\R^N,\|\cdot\|,k}$ is the minimal perimeter of $k$-clusters in $\R^N$ with perimeter given by the norm $\|\cdot\|$, that is,
\[
\J_{\R^N,\|\cdot\|,k}(v) = \inf \Big\{ P_{\|\cdot\|}(\F):\, \F\subseteq \R^N,\, |\F|=v \Big\}\,.
\]
Then, in~\cite{NP2025,FPP} it has been proved that if $M$ is a compact Finsler manifold with a fixed-norm $\|\cdot\|$ (keep in mind that a Riemannian manifold is always a Finsler manifold corresponding to the standard Euclidean norm) then one has
\begin{equation}\label{same-limit}
\lim_{v\to 0}\frac{\J_{M,k}(v)}{\J_{\R^N,\|\cdot\|,k}(v)}= 1\,.
\end{equation}

In view of the above results, the natural question which is left open is whether the connectedness of sufficiently small minimal clusters is true also for Finsler manifolds which are not of fixed-norm type. In this paper we will show that this is not the case, thanks to an explicit example. More precisely, this is our main result.
\maintheorem{main}{Small, non-connected minimal clusters}{For every $N\geq 2$ there exist a smooth, compact, reversible Finsler $N$-dimensional manifold $M$, an integer $k\in\N$, and a sequence of volumes $\eps(j)\in(\R^+)^k$ with $|\eps(j)|\to 0$, such that every $k$-minimal cluster $\E$ with volume $\eps(j)$ is not connected. Moreover, there exists no norm $\|\cdot\|$ on $\R^N$ such that~(\ref{same-limit}) holds.}

We can quickly describe the setting of our example. Our manifold will be a sphere with a smooth Finsler structure which coincides with a flat, Euclidean norm around the north pole, and with (a suitable rescaling of) the $\ell^1$ norm around the south pole; in this setting, we will see that, for a suitable integer $k$ and a suitable large constant $V$, for every sufficiently small $\eps$, any minimal $k$-cluster with volume $(V\eps,\, \eps,\, \dots\, ,\, \eps)$, so with one chamber of volume $V\eps$ and $k-1$ much smaller chambers with volume $\eps$ each, has convenience in having the larger chamber around the north pole and the several smaller ones near the south pole, thus being not connected.\par

The plan of the paper is very simple. In Subsection~\ref{secFinsler} we will briefly present the formal definitions related to Finsler manifolds, and in Subsection~\ref{secCluster} we will describe the isoperimetric problem for clusters fixing all our notation. Then, in Section~\ref{sec:Euclest} we will present the main tools of our construction, which are some estimates for the classical isoperimetric problem for clusters in $\R^N$ with the Euclidean structure. Finally, thanks to these estimates, in Section~\ref{secproof} we will be able to prove Theorem~\mref{main}.

\subsection{Basic definitions for Finsler manifolds\label{secFinsler}}

In this section we will give a short introduction about the geometry of Finsler manifolds. We will only present the tools that we are going to use and to need and we will fix the notation. For a detailed description of Finsler geometry we refer the interested reader to the classical book~\cite{BCS2000}, the more recent book~\cite{Ohtabook}, and the recent papers~\cite{manini2024,ohta2018,ohta2009}, where one can find very clear descriptions, as well as to the introduction of~\cite{NP2025}, where the links between the Finsler geometry and the specific question of the isoperimetric clusters are presented.\par

We can start with the very definition of a smooth Finsler (compact) manifold, which is very classical.

\begin{defin}[Smooth Finsler manifold]\label{defFinMan}
Let $M$ be a compact, connected, smooth $N$-dimensional manifold, and let $\mu$ be a smooth, positive measure on $M$ (this means that $\mu$ can be expressed in local coordinates by means of a strictly positive and smooth density). Moreover, let $j_F: TM \to \R^+$ be a \emph{(smooth) Finsler structure}, that is, a positively $1$-homogeneous function, smooth out of the null section, and such that $j_F^2$ is strictly convex on each tangent space $T_p M$. The triple $(M,j_F,\mu)$ is called a \emph{smooth Finsler manifold}. The manifold is called \emph{reversible} if for any vector $\omega\in TM$ we have $j_F(\omega)=j_F(-\omega)$.
\end{defin}

Having defined what a Finsler manifold is, in order to study isoperimetric problems we need to define what are the volume and the perimeter of a subset $E$ of $M$ (in the next section we will pass to the case of clusters). While the notion of volume is simply given by $\mu$, that is,
\[
|E|_M = \mu(E) = \int_M \Chi{E}(p)\, d\mu(p)\,,
\]
the notion of perimeter requires a bit more care. First of all, the Finsler structure allows to define the length of curves in a natural way. That is, if $\gamma:[0,1]\to M$ is a regular curve, then its length is simply given by
\[
\int_0^1 j_F(\gamma'(t))\,dt\,.
\]
Then, having the length of curves we obtain a (possibly asymmetric) distance on $M$, by setting $d(p,q)$ as the infimum of the lengths of the curves joining $p$ to $q$. Of course, if the Finsler manifold is reversible then the length of a curve does not change upon reverting the orientation, and so $d$ is actually a distance.\par

Having a possibly asymmetric distance, we say that the \emph{Lipschitz constant} of a function $\Phi: M \to \R$ is the smallest constant $L$ such that for every $p,\,q\in M$ one has
\[
\big|\Phi(p)-\Phi(q)\big| \leq L\, \min \big\{ d(p,q),\, d(q,p)\big\}\,,
\]
and the \emph{slope} of a Lipschitz function $\Phi$ at the point $p$ is
\[
|\partial \Phi |(p) = \limsup_{q\to p} \frac{\big(\Phi(q)-\Phi(p)\big)^+}{d(p,q)}\,.
\]
And finally, the \emph{perimeter} of any Borel set $E\subseteq M$ can be defined as
\begin{equation}\label{defP_M}
P_M(E) = \inf \bigg\{ \liminf_{n\to+\infty} \int_M |\partial \Phi_n |(p) \, d\mu(p)\,: \Phi_n \in {\rm Lip}(M),\, \Phi_n \freccia{L^1} \Chi{E}\bigg\}\,.
\end{equation}
It is standard to observe that this notion of perimeter coincides with the standard one whenever $M$ is a Riemannian manifold endowed with its volume measure, and also for a general Finsler manifold it is exactly what one can intuitively think of. Let us be more precise; consider a smooth function $\Phi$ on $M$ and a point $p\in M$, and for every vector $v\in T_p M $ let $\gamma_{p,v}: [0,1]\to M$ be a smooth curve satisfying $\gamma_{p,v}(0)=p$ and $\gamma_{p,v}'(0)=v$. Then, for small $t$ one has
\begin{align*}
d(p,\gamma_{p,v}(t))=t j_F(v)+ o(t)\,, &&\Phi(\gamma_{p,v}(t))-\Phi(p) = t \pr{d\Phi(p)}{v}_M + o(t)\,,
\end{align*}
where with $\pr{\cdot}{\cdot}_M$ we denote the duality between $T_p^* M$ (to which $d\Phi(p)$ belongs) and $T_p M$. As a consequence, by definition the slope of $\Phi$ at $p$ is given by
\[\begin{split}
|\partial \Phi|(p) &= \sup_{v\in T_pM\setminus\{0\}}\lim_{t\to 0}\frac{\left(\Phi(\gamma_{p,v}(t))-\Phi(p)\right)^+}{d(p,\gamma_{p,v}(t))}
=\sup_{v\in T_pM\setminus\{0\}} \frac{\left(\pr{d\Phi(p)}{v}_M\right)^+}{j_F(v)} \\
&=\sup_{v\in T_pM\setminus\{0\}} \Big\langle d\Phi(p), \frac v{j_F(v)}\,\Big\rangle_M = j_F^*(d\Phi(p))\,,
\end{split}\]
where $j_F^* : T^* M \to \R^+$ is dual of $j_F$, defined as usual by
\[
j_F^* (\xi ) = \sup \Big\{ \pr{\xi}{\omega}_M:\, \omega\in T_p M ,\, j_F(\omega)\leq 1\Big\} \qquad \forall \xi \in T^*_p M\,.
\]
As a consequence, we can give an explicit formula to locally describe the perimeter defined in~(\ref{defP_M}). Indeed, let $\varphi: U \to V$ be a smooth chart, hence a diffeomorphism between an open set $U\subseteq \R^N$ and an open set $V \subseteq M$, and let $\rho$ be the density with respect to $\L^N$ of the push-forward measure $(\varphi^{-1})_{\#}\mu$ on $U$. If $E\comp V$ is a set of finite perimeter, then calling $F = \varphi^{-1}(E)\subseteq U$ we can calculate the perimeter of $E$ as
\begin{equation}\label{weuse}
P_M(E)  = \int_{\partial^* F} j_F^*\big(-d\varphi^*_x(\nu_F(x))\big) \rho(x)\, d\H^{N-1}(x)\,,
\end{equation}
where $\nu_F(x)\in\S^{N-1}$ is the outer normal of $F$ at the point $x\in\partial^* F$ (and thus belongs to $T^*_xU$), and $d\varphi^*_x:T^*_xU\to T^*_pV$ is defined as usual as
\[
\pr{d\varphi^* (\omega)}v_V =\pr{\omega}{d\varphi^{-1}(v)}_U
\]
for every $\omega \in T^*_x U$ and $v\in T_pV$.\par

In order to properly investigate isoperimetric clusters in Finsler manifolds there are several other important tools and definitions, as one can see in the introduction of~\cite{NP2025}. However, we will not need anything more in this paper, so we do not introduce them here; just out of completeness, since we have mentioned the ``fixed-norm'' Finsler manifolds above, we can briefly explain what they are. Keeping in mind formula~(\ref{weuse}) above for the perimeter, we can select a particular chart $\varphi$ via a smooth change of variable in such a way that $\varphi^{-1}_\#\mu$ is the Lebesgue measure on $U$. With such a choice of the chart, volume and perimeter of any set $E\comp V$ are given by
\begin{align*}
|E|_M = |F| \,, && P_M(E) = \int_{\partial^* F} h(x,\nu_F(x))\, d\H^{N-1}(x)\,,
\end{align*}
where $|F|$ is the standard Euclidean volume of the set $F\comp U\subseteq\R^N$ and $h:U\times\R^N\to \R^+$ is a suitable function. Thanks to the formula~(\ref{weuse}),
we have an explicit formula for $h$, namely, $h(x,\nu)=  j_F^*\big(-d\varphi^*_x(\nu)\big)$. However, this formula is not particularly interesting by itself. What is important to us, is that $h$ is positively $1$-homogeneous in the second variable, smooth out of the null section $U\times \{0\}$, and the restriction to each section $\{x\}\times \R^N$ has strictly convex sublevels. In other words, writing $\|\omega\| = h(x, \omega)$ for $x=\varphi^{-1}(p)$ and any $\omega\in\R^N$ we have defined a smooth norm $\|\cdot\|$ on $\R^N$. This norm does not only depend on the point $p$ and on the Finsler structure, but also on the choice of the chart $\varphi$; however, the norms corresponding to two different choices of a chart as above coincide up to composition with a unit determinant linear transformation on $\R^N$. As a consequence, we can call such a norm ``infinitesimal norm at $p$'', and it is uniquely defined up to composition with a unit determinant linear map. We can then call a Finsler manifold ``of fixed-norm type'' whenever the infinitesimal norms at all points coincide. In particular, a Riemannian manifold endowed with its Riemannian volume measure is always a fixed-norm Finsler manifold, with all infinitesimal norms coinciding with the standard Euclidean norm (up to unit determinant linear transformations!).

A last comment deserves to be made about the properties of a manifold in the previously known Theorem~\ref{thmknown}, and in our new Theorem~\mref{main}. Since Theorem~\ref{thmknown} shows a property which is true in all the manifolds in a certain class, then the larger the class is the better the theorem is. On the contrary, since our Theorem~\mref{main} shows the existence of a particular manifold in a certain class, this time the theorem becomes stronger if the class is more specific. As a consequence, in the papers~\cite{NP2025,FPP} the Theorem~\ref{thmknown} has been shown in a more general class of manifolds than those of Definition~\ref{defFinMan}, where instead of the smoothness only the ${\rm C}^2$ regularity was required, and the manifolds were not necessarily reversible. Instead, in this paper we stick with the stricter (and more standard) definition where complete smoothness is required, and in addition the example that we build will be a reversible manifold.

\subsection{The isoperimetric problem for clusters\label{secCluster}}

This section is devoted to quickly fix the main notation about the isoperimetric problem for clusters. Given a Finsler manifold $M$ as in Definition~\ref{defFinMan}, for any $k\in\N$ we call \emph{$k$-cluster} any collection of $k$ pairwise disjoint sets in $M$, that we will denote as $\E=(E_1,\, E_2,\, \dots\,,\, E_k)$. The \emph{volume} of a cluster is then the vector of the volumes of the sets $E_i$, so we write
\[
|\E| = \big( |E_1|_M,\, |E_2|_M,\, \dots\, ,\,|E_k|_M \big) \in (\R^+)^k\,.
\]
The \emph{total volume} of the cluster, instead, is the volume of the union $E_0 := \cup_{i=1}^k E_i$, denoted as
\[
\tv\E = \big| E_0 \big|_M = |E_1|_M + |E_2|_M+\ldots +|E_k|_M \in \R^+\,.
\]
Of course, whenever we refer to ``small clusters'', this means clusters with a small total volume.\par

As usual, the \emph{perimeter} of a cluster is defined as
\begin{equation}\label{defperclu}
P_M(\E) = \frac{\sum_{i=0}^k P_M(E_i)} 2\,,
\end{equation}
which is finite if and only if all the chambers (i.e., the sets $E_i$) have finite perimeter. Notice that in the sum there is also the perimeter of the union $E_0$. Notice also the clear meaning of the definition: we call ``perimeter'' of a cluster the measure of all the pieces of boundary, and we divide by $2$ since (almost) every point belongs to two different boundaries; more precisely, an ``internal'' point belongs to the boundary of exactly two different chambers and not of $E_0$, while an ``external'' point belongs to the boundary of a single chamber and also to that of $E_0$.\par

As always, the fact that shared parts of the boundaries are counted just once favours the grouping of the different chambers, and makes it reasonable to think that, in some cases, the cluster is connected (this means that so is the union $E_0$). The goal of this paper is precisely to show a situation when, at any small scale, despite this grouping the union $E_0$ is not connected, differently from what happens in the case of Riemannian or fixed-norm manifolds.

\section{Lower estimates for the Euclidean perimeter of clusters\label{sec:Euclest}}

This section contains the main tools of the paper, that will allow us in the next section to prove Theorem~\mref{main}. We underline that, in this whole section, we will only consider the standard Euclidean perimeter in $\R^N$, so there is no need to know what a Finsler manifold is. To keep the notation light, in this section we will simply write $P$ to denote the standard Euclidean perimeter in $\R^N$, while in the rest of the paper we always write $P_M$ to denote the perimeter in a Finsler manifold $M$, according to the definitions of Subsection~\ref{secFinsler}.\par

The main goal of the section is to give a precise meaning to a simple observation, that is, that a group of many Euclidean balls of unit volume in $\R^N$ cannot be ``perfectly packed'', while many unit cubes can. Let us be more precise; taken any $k$-cluster $\E=(E_1,\, E_2,\, \dots\,,\, E_k)$ in $\R^N$ with volume $(1,\, 1\, \dots\,,\, 1)$, and calling $B$ the (Euclidean) ball with unit volume, then by the definition~(\ref{defperclu}) of perimeter of clusters and by the isoperimetric property of the ball we have the obvious inequality
\begin{equation}\label{trivialineqk}
\p1(\E) \geq \frac 12 \sum_{i=1}^k \p1(E_i) \geq \frac{\p1(B)}2\, k\,.
\end{equation}
The first result of the section will be to show the following improvement of the above inequality.
\begin{lemma}\label{lemma111}
There exists a dimensional constant $c_1=c_1(N)>0$ such that, for every $k\geq 1$ and every $k$-cluster $\E$ with volume $(1,\,1\, \dots\,,\, 1)$, one has
\begin{equation}\label{interthesis}
\p1(\E) \geq \bigg(\frac 12 + c_1\bigg) k\, \p1(B)\,.
\end{equation}
\end{lemma}

A quick comment is in order; since of course the inequality~(\ref{trivialineqk}) is strict for every cluster, in particular a minimal cluster with volume $(1,\, 1\, \dots\, ,\, 1)$, then the validity of~(\ref{interthesis}) with some constant $c_1$ depending on $N$ and on $k$ is obvious. The interesting part of the result is that a constant $c_1$ not depending on $k$ can be found.\par

Let us now consider a more refined case. This time, we let $\E=(F,\, E_1,\, E_2,\, \dots\,, \, E_k)$ be a $(k+1)$-cluster with volume $(V,\, 1,\, 1\, \dots\,,\, 1)$ for some $V$ large enough. Then, again by~(\ref{defperclu}), it is very easy to observe that
\[
\p1(\E) \geq \bigg( V^{\frac{N-1}N}+ \frac k2\bigg)\, \p1(B)\,,
\]
and again the inequality is strict for every $k$ and $V$. Exactly as before, then, we can improve the inequality by a small constant $c_2$, and as before the existence of such a constant depending on $N$, on $k$ and on $V$ is obvious; the interesting part of the next result is that we can provide a constant depending only on $N$, and not on $k$ or on $V$.

\begin{prop}\label{donow}
There exists a dimensional constant $c_2=c_2(N)>0$ such that, for any $V\geq 1$, if $k\geq 1$ and $\E$ is a $(k+1)$-cluster with volume $(V,\,1,\, \dots \,,\, 1)$ then
\begin{equation}\label{extraperimeter}
\p1(\E) \geq \bigg( V^{\frac{N-1}N}+ \bigg( \frac 12 + c_2\bigg) k \bigg)\, \p1(B)\,.
\end{equation}
\end{prop}

The rest of this section is divided in four parts. In the first one we present the statement of an important (and quite well-known) quantitative form of the isoperimetric inequality, due to Fusco and Julin, that we will need later on. Then, in the second part, we present an immediate geometrical fact, and we deduce an estimate on the common boundary between two sets. Thanks to this estimate, in the last two parts we will be able to show first Lemma~\ref{lemma111} and then Proposition~\ref{donow}.

\subsection{A quantitative result}

In this section we present the statement of a celebrated strong form of the quantitative isoperimetric inequality, due to Fusco and Julin (see~\cite{FJ}). Let $E\subseteq \R^N$ be a set with finite perimeter and with $0<|E|<+\infty$, and let $B_E$ denote the ball centered in the origin and with the same volume as $E$, hence so that $P(E)\geq P(B_E)$, with strict inequality unless $E$ coincides with a translation $y+B_E$ of $B_E$. We call \emph{isoperimetric deficit} $\delta(E)$ and \emph{asymmetry} $\lambda(E)$ the quantities
\begin{equation}\label{eq:asymmetry}
\begin{array}{c}
\delta(E) :=\bal \frac{\p1(E)-\p1(B_E)}{\p1(B_E)}\,,\eal\\
\lambda(E) := \inf\limits_{y\in\R^N}\bal\frac{|E\Delta (y+B_E)|}{|B_E|} +\Bigg(\frac 1{|B_E|^{\frac{N-1}N}}\int_{\partial^*E} \big|\nu_E(x)-\theta_{yx}\big|^2d\H^{N-1}(x)\Bigg)^{1/2}\eal\,.
\end{array}
\end{equation}
The meaning of the isoperimetric deficit is quite obvious, as well as the fact that it is always positive, and strictly positive unless $E$ is a translation of $B_E$. Concerning the asymmetry, here and in the rest of the paper whenever $x\neq y$ we denote by $\theta_{yx}$ the direction of the vector pointing from $y$ to $x$, namely,
\[
\theta_{yx}=\frac{x-y}{|x-y|} \in \S^{N-1}\,.
\]
Also the meaning of the asymmetry is then easy to understand: the first part is the volume of the symmetric difference between $E$ and the ``best placed'' translation $y+B_E$ of $B_E$, hence it is again always positive, and strictly positive unless $E$ actually coincides with $y+B_E$. The second part, instead, is the average ``deviation'' between the outer normal vector $\nu_E(x)$ at a point $x$ of the boundary $\partial^* E$, and the ``expected normal'' $\theta_{yx}$; observe that also this second part is always positive, and it is strictly positive unless $E$ is a ball centered at $y$.\par

It is standard to notice that the infimum in the definition of the asymmetry is in fact a minimum. We will then call \emph{Fraenkel point for $E$} any point $y$ realizing the minimum in the definition of the asymmetry. The quantitative isoperimetric inequality reads then as follows.

\begin{theorem}[\cite{FJ}, Theorem~1.1]\label{sqii}
For any $N\geq 2$ there is a constant $c_3 = c_3(N)>0$ such that for every set $E\subseteq \R^N$ of finite perimeter one has
\begin{equation}\label{eq:fusco-julin-quantitative}
\delta(E) \geq c_3 \lambda(E)^2\,.
\end{equation}
\end{theorem}

\subsection{A first geometrical result}

We start with the following simple observation.
\begin{lemma}\label{basiclemma}
Let $\sigma\in(0,1]$ and $R>0$. If $z$ and $w$ are two points in $\R^N$, and $x$ is a point having distance at least $\sigma R$ from a line passing through $z$ and $w$, and such that $|x-z|\leq R$, then one has $\big|\theta_{zx}+\theta_{wx}\big|\geq \sigma$.
\end{lemma}
\begin{proof}
First of all, we notice that the case when $z=w$ is trivial, because in this case $\theta_{zx}=\theta_{wx}$ for any $x$, and then $|\theta_{zx}+\theta_{wx}|=2$. Suppose then that $z\neq w$, and let $x$ be any point as in the claim. Up to a translation and/or a rotation, we can assume that $z=(0,\,0,\, \dots\,,\, 0)$, $w=(a,\, 0,\, 0,\, \dots\,,\, 0)$ and $x=(x_1,\, x_2,\, 0,\, 0,\, \dots\,,\, 0)$, with some $x_1\geq 0$ and some $x_2\geq \sigma R$. We then have that $|\theta_{zx}+\theta_{wx}|\geq |\theta_{zx}-{\rm e}_1|$, and the latter term is minimal if $x_2=\sigma R$ and $|x|=R$. And finally, in this particular case we get
\[
|\theta_{zx}-{\rm e}_1| = 2 \sin \bigg( \frac{\arcsin (\sigma)}2\bigg) \geq 2\sin \bigg( \frac{\arcsin (\sigma)}2\bigg)\cos \bigg( \frac{\arcsin (\sigma)}2\bigg) = \sin \left( \arcsin (\sigma)\right) = \sigma\,.
\]
\end{proof}

As a consequence of the above simple result, we can show the following estimate.

\begin{lemma}\label{secondbasic}
For every sufficiently small $\sigma>0$ there exists $c_4=c_4(N,\sigma)>0$ so that the following holds. Let $F,\, E$ be two disjoint sets, with $|E|=1$, and such that
\[
\H^{N-1}\big(\partial^* E\cap\partial^* F\big) \geq \sigma \p1(B)\,,
\]
and let $w\in\R^N$ be any given point. Then, at least one of the two following estimates is true
\begin{gather}
\p1(E) \geq (1+c_4) \p1(B)\,, \label{caseA}\\
\int_{\partial^* F\cap\partial^* E} \Big| \nu_F(x) - \theta_{wx}\Big|^2\, d\H^{N-1}(x)\geq \frac{c_4}{c_3}\,. \label{caseB}
\end{gather}
In particular, thanks to Theorem~\ref{sqii}, if also $|F|=1$ then $\max\{\p1(E),\, \p1(F) \} \geq (1+c_4) \p1(B)$.
\end{lemma}
\begin{proof}
For simplicity of notations we can assume that the origin is a Fraenkel point for $E$, that is, a point realizing the minimum in the definition~(\ref{eq:asymmetry}) of the asymmetry. Then, we have to distinguish a few of possible different cases. Keep in mind that we are calling $B$ the Euclidean ball centered at the origin and with unit volume, hence with radius $\omega_N^{-1/N}$; then, we call now $B_2$ the ball, still centered in the origin, whose radius is twice as large, hence $2\omega_N^{-1/N}$
\case{I}{If $E$ has much perimeter outside of $B_2$.}
First of all, let us assume that
\[
\p1 \big(E ; \R^N\setminus B_2\,\big) \geq \frac\sigma 2 \p1(B)\,.
\]
We claim that there exists a constant $c_4'=c_4'(N,\sigma)$, only depending on $N$ and $\sigma$, such that~(\ref{caseA}) holds true for $E$ under this assumption with $c_4'$ in place of $c_4$. Indeed, suppose that this is not true; then, there exists a sequence $\{ E_n\}_{n\in\N}$ of sets such that $\p1(E_n)\to \p1(B)$, each $E_n$ has volume $1$ and admits the origin as a Fraenkel point, and $\p1(E_n; \R^N\setminus B_2) \geq \sigma \p1(B)/2$. The fact that $\p1(E_n)\to \p1(B)$ implies that $|E_n \Delta B|\to 0$ by Theorem~\ref{sqii}, and then by lower semi-continuity of the perimeter we get
\[\begin{split}
\p1(B) &=\p1(B; B_2) \leq \liminf \p1(E_n; B_2) = \liminf \Big( \p1(E_n) - \p1(E_n; \R^N\setminus B_2)\Big)\\
&= \p1(B) -\limsup \p1(E_n; \R^N\setminus B_2)
\leq \p1(B) \bigg( 1 -\frac \sigma 2\bigg)\,,
\end{split}\]
which is impossible. The contradiction shows then the existence of the constant $c_4'$ as in our claim.

Notice that, if we are not in Case~I, then by our standing hypothesis $\H^{N-1}\big(\partial^* E\cap\partial^* F\big) \geq \sigma \p1(B)$ we get
\[
\H^{N-1}\big( \partial^* E \cap \partial^* F \cap B_2\big) \geq \frac \sigma 2\, \p1(B)\,.
\]

\case{II}{If there is some line around which too much perimeter concentrates.}
Let us now assume that there exists a line $\ell$ passing through a Fraenkel point for $E$ such that
\[
\H^{N-1}\Big( \partial^* E\cap \C_\ell\Big) > \frac \sigma 4 \, \p1(B)\,,
\]
where $\C_\ell$ is the set of points in $\R^N$ having distance less than $\sigma\omega_N^{-1/N}/6$ from the line $\ell$. Then, we claim that there exists a constant $c_4''=c_4''(N,\sigma)$, only depending on $N$ and on $\sigma$, such that again~(\ref{caseA}) holds true for $E$ under this assumption, with $c_4''$ in place of $c_4$. The proof is very similar to the one of the previous case; namely, otherwise there exists a sequence $\{E_n\}$ of sets with $\p1(E_n)\to \p1(B)$, each of them having unit volume and admitting the origin as a Fraenkel point, and with $\H^{N-1}(\partial^* E_n \cap \C_\ell)> \sigma \p1(B)/4$ for some fixed $\C_\ell$. We deduce that $|E_n\Delta B|\to 0$, and by lower semicontinuity we get 
\[\begin{split}
\p1(B;\overline{\C_\ell})&= \p1(B) - \p1(B;\R^N\setminus\overline{\C_\ell}) \geq \p1(B) - \liminf \p1(E_n;\R^N\setminus\overline{\C_\ell})\\
&= \p1(B) - \liminf \Big( \p1(E_n) - \p1(E_n;\overline{\C_\ell})\Big) = \limsup \p1(E_n;\overline{\C_\ell}) \geq \frac\sigma 4\, \p1(B)\,.
\end{split}\]
However, if $\sigma$ is small then $\p1(B;\overline{\C_\ell}) \approx 2 \omega_{N-1} (\sigma\omega_N^{-1/N}/6)^{N-1}$, and this gives the desired contradiction with the above inequality. Indeed, if $N=2$ the contradiction for a sufficiently small $\sigma$ comes from the fact that $4/6<\pi/2$, while if $N\geq 3$ it comes since $\sigma^{N-1}\ll \sigma$. Summarizing, the existence of the constant $c_4''$ as in our claim is proved.

Now, let $w$ be any point in $\R^N$, let $\ell$ be the line passing through $w$ and the origin (or any such line if $w$ is the origin), and let us call
\[
\Gamma = \partial^* E \cap \partial^* F \cap B_2 \setminus \C_\ell\,.
\]
If we are neither in Case~I nor in Case~II, then
\begin{equation}\label{not12}
\H^{N-1}(\Gamma ) \geq \frac \sigma 4\, \p1(B)\,.
\end{equation}
Notice that, since $\Gamma\subset B_2\setminus \C_{\ell}$, we can apply Lemma~\ref{basiclemma} with $z=0$ and $R=2\omega_N^{-1/N}$, any point $x\in\Gamma$, and with $\sigma/12$ in place of $\sigma$, and this tells us that
\[
|\theta_{zx}+\theta_{wx}|\geq \frac \sigma{12}\,.
\]
Since $x\in \partial^* E \cap \partial^* F$ and the two sets $E$ and $F$ are disjoint, we have that $\nu_E(x) = -\nu_F(x)$; as a consequence, we get
\[
\big| \nu_E(x) - \theta_{zx} \big| + \big|\nu_F(x) -\theta_{wx} \big| \geq \frac \sigma{12}\,,
\]
which implies that at least one between $|\nu_E(x)-\theta_{zx}|$ and $|\nu_F(x)-\theta_{wx}|$ is larger than $\sigma/24$. Thus, calling
\begin{align*}
\Gamma_1 = \bigg\{ x \in \Gamma,\, |\nu_E(x)-\theta_{zx}|\geq \frac \sigma{24}\bigg\}\,, &&
\Gamma_2 = \bigg\{ x \in \Gamma,\, |\nu_F(x)-\theta_{wx}|\geq \frac \sigma{24}\bigg\}\,,
\end{align*}
by~(\ref{not12}) we deduce
\[
\H^{N-1}(\Gamma_1)+\H^{N-1}(\Gamma_2) \geq \frac \sigma 4\, \p1(B)\,.
\]
We can then finally divide the possible situations in two further subcases.

\case{III}{If $\H^{N-1}(\Gamma_1)\geq \sigma \p1(B)/8$.}
Let us first suppose that $\H^{N-1}(\Gamma_1)$ is larger than $\sigma \p1(B)/8$. In this case, since $z=0$ is a Fraenkel point for $E$, the quantitative isoperimetric Theorem~\ref{sqii} together with the definition~(\ref{eq:asymmetry}) of the asymmetry gives
\[
\delta(E) \geq c_3 \lambda(E)^2 \geq c_3 \int_{\partial^*E} \big|\nu_E(x)-\theta_{zx}\big|^2d\H^{N-1}(x)
\geq c_3 \cdot\,\frac{\sigma^3 \p1(B)}{8\cdot 24^2}\,,
\]
and then also in this case~(\ref{caseA}) holds true with $c_4'''$ in place of $c_4$, where
\[
c_4'''=c_4'''(N,\sigma) := c_3 \cdot\,\frac{\sigma^3 \p1(B)}{8\cdot 24^2}\,.
\]

\case{IV}{If $\H^{N-1}(\Gamma_2)\geq \sigma \p1(B)/8$.}
The last case to be considered is when $\H^{N-1}(\Gamma_2)$ is larger than $\sigma \p1(B)/8$. And in this case, by construction we have
\[
\int_{\partial^* F\cap\partial^* E} \Big| \nu_F(x) - \theta_{wx}\Big|^2\, d\H^{N-1}(x)\geq
\int_{\Gamma_2} \Big| \nu_F(x) - \theta_{wx}\Big|^2\, d\H^{N-1}(x)\geq
\frac{\sigma^3 \p1(B)}{8\cdot 24^2} = \frac{c_4'''}{c_3}\,.
\]
Hence, this time~(\ref{caseB}) holds true with $c_4'''$ in place of $c_4$. Summarizing, the first part of the statement is established up to defining $c_4= \min \big\{c_4',\, c_4'',\, c_4'''\big\}$.\par

In particular, if also $|F|=1$, then applying the first part of the statement with $w$ being a Fraenkel point of $F$ we obtain the thesis directly by Theorem~\ref{sqii}.
\end{proof}

\subsection{The proof of Lemma~\ref{lemma111}}

This section is devoted to present the proof of our first estimate on the Euclidean perimeter of clusters. In the following, we will denote as usual $B(x,r)$ the ball centered at $x\in\R^N$ and with radius $r>0$, so that the unit ball is $B=B(0,\omega_N^{-1/N})$.
\proofof{Lemma~\ref{lemma111}}
Let $\E=(E_1,\, E_2,\,\dots\,,\, E_k)$ be any cluster with volume~$(1,\,1\, \dots\,,\, 1)$. For every $1\leq i \leq k$, we let $x_i$ be a Fraenkel point for $E_i$, and we denote by $B^i=B(x_i,\omega_N^{-1/N})$ the corresponding ball. Moreover, we call $B^i_+=B(x_i,2\omega_N^{-1/N})$ the ball with the same center but a radius twice bigger; moreover, we write
\[
\widetilde\partial^* E_i = \partial^* E_i \cap B_+^i\,.
\]
Arguing exactly as in Case~I of Lemma~\ref{secondbasic}, there exists a constant $c_5=c_5(N)$ such that
\begin{equation}\label{constc'}
\H^{N-1}\big(\widetilde\partial^* E_i\big) < \frac{24}{25}\, \p1(B) \quad \Longrightarrow \quad \p1(E_i)\geq \p1(B) (1 +c_5)\,.
\end{equation}
We can then define
\[
I = \Big\{ 1\leq i\leq k,\, \p1(E_i)< \p1(B) (1+ c_5)\Big\}\,.
\]
Let us consider $i\in I$; by~(\ref{eq:fusco-julin-quantitative}) we have
\[
\big| E_i \Delta B^i| \leq \sqrt {\frac{c_5}{c_3}}\,.
\]
In particular, if $i,\, j\in I$ then we have
\begin{equation}\label{ell}
\big|B^i\cap B^j\big| \leq \sqrt {\frac{c_5}{c_3}}\,, 
\end{equation}
since all the points of $B^i\cap B^j$ belong to $B^i\setminus E_i$ or to $B^j\setminus E_j$. Now, for any index $i\in I$ we define
\[
J_i = \Big\{j \in I\setminus \{i\},\, B_+^i \cap B_+^j \neq \emptyset\Big\}\,.
\]
By construction, for every $j\in J_i$ the ball $B^j$ is entirely contained in the ball $B(x_i,5\omega_N^{-1/N})$ with the same center as $B^i$ and radius $5$ times bigger. In addition, for any two different indices $j,\, j'\in J_i$ we can apply~(\ref{ell}) to get that the intersection $B^j\cap B^{j'}$ is small. Therefore, there exists a purely dimensional number $C=C(N)$ such that
\[
\# J_i \leq C \qquad \forall\, i\in I\,.
\]
We now subdivide the indices of $I$ in four subgroups as follows
\begin{gather*}
I_1 = \bigg\{ i \in I,\, \H^{N-1}\Big(\widetilde\partial^* E_i\cap \bigcup\nolimits_{j\in I\setminus \{i\}} \widetilde\partial^* E_j\Big) \geq \frac 6{25}\,\p1(B)\bigg\}\,,\\
I_2 = \bigg\{ i \in I\setminus I_1,\, \H^{N-1}\Big(\widetilde\partial^* E_i\cap \bigcup\nolimits_{j\in I\setminus \{i\}} \partial^* E_j\setminus \widetilde\partial^* E_j\Big) \geq \frac 6{25}\,\p1(B)\bigg\}\,,\\
I_3 = \bigg\{ i \in I\setminus (I_1\cup I_2),\, \H^{N-1}\Big(\widetilde\partial^* E_i\cap \bigcup\nolimits_{j\notin I} \partial^* E_j\Big) \geq \frac 6{25}\,\p1(B)\bigg\}\,,\\
I_4 = I \setminus \big( I_1\cup I_2\cup I_3)\,.
\end{gather*}
Now, let us fix any $i\in I_1$; first of all, we notice that it must be
\[
\H^{N-1}\Big(\widetilde\partial^* E_i\cap \bigcup\nolimits_{j\in J_i} \widetilde\partial^* E_j\Big) \geq \frac 6{25}\,\p1(B)\,,
\]
because by definition if $j\neq i$ is an index in $I\setminus J_i$ then $B^i_+\cap B^j_+=\emptyset$, and then $\widetilde\partial^* E_i\cap\widetilde \partial^* E_j=\emptyset$. Then, since $J_i$ contains at most $C$ elements, there exists some $j\in J_i$ such that
\[
\H^{N-1}\big(\widetilde\partial^* E_i\cap \widetilde\partial^* E_j\big) \geq \frac 6{25C}\, \p1(B)\,.
\]
Let us now select a constant $\sigma$, only depending on $N$, small enough so that Lemma~\ref{secondbasic} can be applied and such that $\sigma\leq 6/25C$. Since both $E_i$ and $E_j$ have unit volume, we apply Lemma~\ref{secondbasic} to $E_i$ and $E_j$, and we deduce that
\begin{equation}\label{unodeidue}
\max \Big\{ \p1(E_i),\, \p1(E_j)\Big\} \geq (1+c_4) \p1(B)\,.
\end{equation}
Notice that the constant $c_4$ depends only on $N$ and on $\sigma$, thus actually only on $N$ since $\sigma$ has been chosen only depending on $N$.\par

We are now going to prove the validity of~(\ref{interthesis}) by distinguishing five possible cases. More precisely, we call for brevity $I_0 = \{1,\,2,\,\dots\,,\, k\}\setminus I$, so that $\{1,\,2,\,\dots\,,\,k\}$ is the disjoint union of $I_0,\, I_1,\, I_2,\, I_3$ and $I_4$; then, either $I_0$ has at least $k/25$ elements, or one between $I_1,\, I_2,\, I_3,\, I_4$ has at least $6k/25$ elements. We consider separately the five possibilities.

\case{I}{If $\# I_0\geq k/25$.}
This is the easiest case to deal with. In fact, for every $i\in I_0$ we have $\p1(E_i)\geq \p1(B)(1 + c_5)$ by definition, and then
\[\begin{split}
\p1(\E) &> \frac 12\, \sum_{i=1}^k \p1(E_i) = \frac 12\, \bigg( \sum_{i\in I_0} \p1(E_i) + \sum_{i\notin I_0} \p1(E_i)\bigg)\\
&\geq \frac 12 \bigg( \# I_0 \big(\p1(B)(1+c_5) \big) + (k-\# I_0) \p1(B)\bigg)
= \frac{\p1(B)}2\, \bigg( k + c_5\# I_0\bigg)\\
&\geq \bigg( \frac 12 +\frac{c_5}{50}\bigg) k \p1(B)\,.
\end{split}\]
The validity of~(\ref{interthesis}) is then established in this first case, as soon as $c_1\leq c_1':= c_5/50$.

\case{II}{If $\# I_1\geq 6k/25$.}
Let us now consider this second case. As observed above, for every $i\in I_1$ there exists some $j\in J_i$ such that~(\ref{unodeidue}) holds. Keeping in mind that, for $i,\, j \in I$, one has that $j\in J_i$ if and only if $i\in J_j$, this immediately implies that
\[
\# \widetilde I \geq \frac{6k}{25(C+1)}\,, \qquad \hbox{where} \qquad \widetilde I =\bigg\{ i \in I,\, \p1(E_i) \geq (1+c_4) \p1(B)\bigg\}\,.
\]
And then, in this case we get
\[\begin{split}
\p1(\E) &> \frac 12\, \sum_{i=1}^k \p1(E_i) = \frac 12\, \bigg( \sum_{i\in \widetilde I} \p1(E_i) + \sum_{i\notin \widetilde I} \p1(E_i)\bigg)\\
&\geq \frac 12\, \bigg(\# \widetilde I \big(1+c_4\big) \p1(B) + (k-\#\widetilde I) \p1(B)\bigg)
= \frac {\p1(B)}2\, \Big(k + c_4 \# \widetilde I\Big)\\
&\geq \bigg( \frac 12+\frac{3 c_4}{25(C+1)} \bigg)\, k \p1(B)\,,
\end{split}\]
so once again the validity of~(\ref{interthesis}) follows as soon as $c_1\leq c_1'':=3 c_4/(25(C+1))$.

\case{III}{If $\# I_2\geq 6k/25$.}
We pass now to the third case. We have
\[
\bigcup_{j\in I} \big(\partial^* E_j \setminus \widetilde\partial^* E_j\big)\supseteq \bigcup_{i\in I_2} \Big(\widetilde\partial^* E_i \cap \bigcup\nolimits_{j\in I\setminus \{i\}} \partial E_j\setminus \widetilde\partial^* E_j\Big)\,,
\]
and it follows that
\[
\sum_{j\in I} \H^{N-1}\Big(\partial^* E_j \setminus \widetilde\partial^* E_j\Big) \geq \# I_2\, \frac 6{25}\,\p1(B)> \frac 6{125}\, k\p1(B)\,.
\]
However, for every $j\in I$ by~(\ref{constc'}) we know that $\H^{N-1}(\widetilde\partial^* E_j)\geq 24 \p1(B)/25$, and then
\[\begin{split}
\sum_{j=1}^k \p1(E_j)&= \sum_{j\in I_0} \p1(E_j) + \sum_{j\in I} \p1(E_j)\\
&= \sum_{j\in I_0} \p1(E_j)+\sum_{j\in I} \H^{N-1}\big(\widetilde\partial^* E_j\big)+ \sum_{j\in I}\H^{N-1}\big(\partial^* E_j\setminus\widetilde\partial^* E_j\big)\\
&\geq \# I_0 \p1(B) + \# I \, \frac{24}{25}\, \p1(B) + \frac 6{125}\, k \p1(B)
\geq \bigg(\frac {24}{25} + \frac 6{125}\bigg) k\p1(B)
=\frac{126}{125}\, k \p1(B)\,.
\end{split}\]
Thus, also in this case we get~(\ref{interthesis}) as soon as $c_1\leq c_1''':= 1/250$, since
\[
\p1(\E)>\frac 12\, \sum_{i=1}^k \p1(E_i) \geq \bigg(\frac 12 + \frac 1{250}\bigg) k\p1(B)\,.
\]

\case{IV}{If $\# I_3\geq 6k/25$.}
We consider now the fourth case. In this case, we have
\[
\sum_{i\notin I} \H^{N-1}(\partial^* E_i) \geq \# I_3 \,\frac 6 {25}\, \p1(B)> \frac 6 {125}\, k \p1(B)\,.
\]
We can now assume that $\# I_0< k/25$ (and then, $\# I\geq 24 k /25$), because otherwise the thesis has been already obtained in Case~I. As a consequence, this time we have
\[\begin{split}
\p1(\E) &> \frac 12\, \sum_{i=1}^k \p1(E_i) = \frac 12 \bigg(\sum_{i\notin I} \p1(E_i) + \sum_{i\in I} \p1(E_i)\bigg)
> \frac 12 \bigg( \frac 6 {125}\, k \p1(B) + \# I \,\p1(B)\bigg)\\
&\geq\frac 12 \bigg( \frac 6 {125}\, k \p1(B) + \frac{24}{25}\, k \p1(B)\bigg)
= \bigg( \frac 12 + \frac 1{250}\bigg) k \p1(B)\,,
\end{split}\]
and so again~(\ref{interthesis}) follows if $c_1\leq c_1'''$.

\case{V}{If $\# I_4\geq 6k/25$.}
Let us now finally consider the last case. Since for every $i\in I$ we have $\H^{N-1}(\widetilde\partial^* E_i) \geq 24 \p1(B)/25$, if $i\in I_4$ this means that
\[
\H^{N-1}\Big(\widetilde\partial^* E_i\setminus \bigcup_{j\neq i} \partial^* E_j\Big) \geq \frac 6{25}\, \p1(B)\,.
\]
However, we have that
\[
\widetilde\partial^* E_i\setminus \bigcup_{j\neq i} \partial^* E_j \subseteq \partial^* E_0\,,
\]
where $E_0=\cup_{j=1}^k E_j$. As a consequence, this time we evaluate
\[
\p1(\E)=\frac 12\, \bigg(\p1(E_0) + \sum_{i=1}^k \p1(E_i)\bigg)
\geq \frac 12\, \bigg(\# I_4\, \frac 6{25}\, \p1(B)+ k \p1(B)\bigg)
> \bigg(\frac 12+ \frac 3{125} \bigg)\,k \p1(B)\,,
\]
and so once again we get the validity of~(\ref{interthesis}) if $c_1\leq c_1'''':=3/125$. Summarizing, we have completed the proof of the validity of~(\ref{interthesis}) defining $c_1=\min\{c_1',\, c_1'',\, c_1''',\, c_1''''\}$, and as noticed at the beginning this is equivalent to the thesis.
\end{proof}

\subsection{The proof of Proposition~\ref{donow}}
This section is devoted to present the proof of our second estimate on the Euclidean perimeter of the clusters, for the case when the cluster has $k$ chambers with equal volume, plus another, much bigger, chamber.

\proofof{Proposition~\ref{donow}}
We consider a given $(k+1)$-cluster $\E=(F,E_1,\, E_2,\, \dots\,,\, E_k)$ with volume $(V,\,1,\, \dots \,,\, 1)$. For brevity of notations, we call $\widetilde \E$ the $k$-cluster made by the $k$ sets $E_i$, and $E^+=\cup_{i=1}^kE_i$ the union of these sets. In addition, we call
\begin{align*}
\Gamma^- = \partial^* E^+ \cap \partial^* F\,, && \Gamma^+ = \partial^* E^+ \setminus \partial^* F\,.
\end{align*}
We can then prove the claim considering separately a few possibilities.

\case{I}{If $\Gamma^-$ is small.}
Let us first assume that
\[
\H^{N-1}(\Gamma^-) \leq \frac{c_1}2\, k \p1(B)\,,
\]
where $c_1=c_1(N)$ is the constant of Lemma~\ref{lemma111}. In this case, applying Lemma~\ref{lemma111} to the cluster $\widetilde \E$ we have
\[
\p1(\E) = \p1(F) + \p1(\widetilde \E) - \H^{N-1}(\Gamma^-)
\geq \bigg( V^{\frac{N-1}N} + \bigg(\frac 12 + \frac{c_1}2 \bigg) \, k \bigg)\p1(B)\,,
\]
so in this first case~(\ref{extraperimeter}) is established if $c_2\leq c_2' : =c_1/2$.

\case{II}{If $\Gamma^-$ is larger than $\Gamma^+$.}
We now assume that for some small $\xi$, only depending on $N$ and to be specified later, one has
\[
\H^{N-1}(\Gamma^-) \geq \H^{N-1}(\Gamma^+) + 2\xi\, k \p1(B)\,.
\]
In this case, directly by definition we get
\[\begin{split}
\p1(\E) &= \p1(F\cup E^+) + \frac 12\, \sum_{i=1}^k \p1(E_i) + \frac 12\, \Big( \H^{N-1}(\Gamma^-)-\H^{N-1}(\Gamma^+)\Big)\\
&\geq \bigg( V^{\frac{N-1}N} + \bigg(\frac 12 +\xi \bigg) k\bigg) \, \p1(B)\,,
\end{split}\]
and then also in this second case~(\ref{extraperimeter}) is established as soon as $c_2\leq c_2'':= \xi$.

Notice that, if we are neither in Case~I nor in Case~II, then
\begin{equation}\label{neitherInorII}
\frac{c_1}2\, k \p1(B)<\H^{N-1}(\Gamma^-) < \H^{N-1}(\Gamma^+) + 2\xi\, k \p1(B)\,.
\end{equation}
Moreover, we can simply estimate
\begin{equation}\label{umt}\begin{split}
\p1(\E) &= \p1(F) + \frac 12\, \sum_{i=1}^k P(E_i) + \frac 12 \Big( \H^{N-1}(\Gamma^+)-\H^{N-1}(\Gamma^-)\Big)\\
&\geq \p1(F) + \frac 12\, \sum_{i=1}^k P(E_i) - \xi\, k \p1(B)\,.
\end{split}\end{equation}
We now define
\[
I = \bigg\{ 1\leq i \leq k,\, \H^{N-1}(\partial^* E_i\cap \partial^* F)\geq \frac {c_1}6\, \p1(B)\bigg\}\,,
\]
and we further subdivide the possible situations in more cases.

\case{III}{If $I$ has few elements.}
Let us first assume that
\begin{equation}\label{assfew}
\# I \leq \frac{c_1}4\, k\,.
\end{equation}
In this case, we observe that
\[
\sum_{i\notin I} \H^{N-1}(\partial^* E_i \cap \partial^* F)\leq \frac {c_1}6\, k \p1(B)\,,
\]
and then the validity of~(\ref{neitherInorII}) ensures that
\[
\sum_{i\in I} \p1(E_i) \geq \sum_{i\in I} \H^{N-1}(\partial^* E_i \cap \partial^* F) = \H^{N-1}(\Gamma^-) - \sum_{i\notin I} \H^{N-1}(\partial^* E_i \cap \partial^* F)
\geq \frac {c_1}3\, k \p1(B)\,.
\]
Inserting this estimate into~(\ref{umt}) and using the assumption~(\ref{assfew}), we get
\[\begin{split}
\p1(\E) &\geq \p1(F) + \frac 12\, \sum_{i\in I} P(E_i) + \frac 12\, \sum_{i\notin I} P(E_i) - \xi\, k \p1(B)\\
&\geq \bigg( V^{\frac{N-1}N} + \frac{c_1}6\, k + \frac {k-\#I}2 -\xi\,k \bigg) \,\p1(B)
\geq \bigg( V^{\frac{N-1}N} + \bigg( \frac 12 + \frac {c_1}{24}-\xi\bigg)\, k\bigg)\, \p1(B)\,.
\end{split}\]
Therefore, up to choosing $\xi<c_1/24$, also this time we get~(\ref{extraperimeter}) provided $c_2\leq c_2''':= c_1/24 - \xi$.

We are now left to consider the case when
\begin{equation}\label{assenough}
\# I > \frac{c_1}4\, k\,.
\end{equation}
Now, we select a constant $\sigma$, only depending on $N$, small enough so that Lemma~\ref{secondbasic} can be applied and such that $\sigma\leq c_1/6$. With this choice of $\sigma$, the value of the constant $c_4$ of Lemma~\ref{secondbasic} depends only on $N$. In addition, we fix a Fraenkel point $w$ for $F$.

Let $i\in I$ be any index. By definition of $I$ and since $\sigma\leq c_1/6$ we can apply Lemma~\ref{secondbasic} to the sets $F$ and $E=E_i$, and we get that either $\p1(E_i)\geq (1+c_4)\p1(B)$, or
\begin{equation}\label{laleth}
\int_{\partial^* F\cap\partial^* E_i} \Big| \nu_F(x) - \theta_{wx}\Big|^2\, d\H^{N-1}(x)\geq \frac{c_4}{c_3}\,.
\end{equation}
We can then subdivide $I$ into two parts, namely,
\begin{align*}
I' = \big\{ i\in I,\, \p1(E_i)\geq (1+c_4) \p1(B)\big\}\,, && I'' = I\setminus I'\,.
\end{align*}
As just observed, we know that~(\ref{laleth}) is true for every $i\in I''$. We conclude the proof by separately considering what happens when either $I'$ or $I''$ has more elements.

\case{IV}{If $I'$ has many elements.}
We now consider the situation if
\begin{equation}\label{assenough'}
\# I' > \frac{c_1}8\, k\,.
\end{equation}
In this case, thanks to~(\ref{umt}) we have
\[\begin{split}
\p1(\E)&\geq \p1(F) + \frac 12\, \sum_{i=1}^k P(E_i) - \xi\, k \p1(B)
= \p1(F) + \frac 12\, \sum_{i\in I'} P(E_i)+\frac 12\, \sum_{i\notin I'} P(E_i) - \xi\, k \p1(B)\\
&\geq \bigg(V^{\frac{N-1}N} + \frac{1+c_4}2\, \#I'+ \frac{k - \# I'}2 -\xi\, k\bigg) \,\p1(B)
\geq \bigg(V^{\frac{N-1}N} + \frac k 2+ \frac{c_4 c_1}{16}\, k -\xi\, k\bigg) \,\p1(B).
\end{split}\]
Thus, up to choosing $\xi< c_4 c_1/16$, once again we get the validity of the inequality~(\ref{extraperimeter}), provided that $c_2\leq c_2'''':=c_4 c_1/16 - \xi$.

\step{V}{If $I''$ has many elements.}
Since we have reduced ourselves to consider the case when~(\ref{assenough}) holds, and we have already discussed what happens if~(\ref{assenough'}) holds, the last situation to be considered is if
\[
\# I'' > \frac{c_1}8\, k\,.
\]
Since for every $i\in I''$ we have the validity of~(\ref{laleth}), and the different sets $\partial^* F \cap \partial^* E_i$ are pairwise essentially disjoint, in this situation we have
\[
\int_{\partial^* F} \Big| \nu_F(x) - \theta_{wx}\Big|^2\, d\H^{N-1}(x) \geq \frac{c_1 c_4}{8c_3}\, k\,.
\]
By Theorem~\ref{sqii} we have then
\[\begin{split}
\p1(F) &= \p1(B) V^{\frac{N-1}N} \big( 1 +\delta(F)\big)
\geq \p1(B) V^{\frac{N-1}N} \big( 1 + c_3 \lambda(F)^2\big)\\
&\geq \p1(B) V^{\frac{N-1}N} \bigg( 1 + \frac {c_3}{V^{\frac{N-1}N}}\int_{\partial^*F} \big|\nu_F(x)-\theta_{wx}\big|^2d\H^{N-1}(x)\bigg)\\
&\geq \p1(B) V^{\frac{N-1}N} \bigg( 1 + \frac {c_1 c_4 }{8V^{\frac{N-1}N}}\, k\bigg)
= \p1(B) V^{\frac{N-1}N} + \frac {c_1 c_4}{8}\, k \p1(B)\,.
\end{split}\]
Inserting this estimate in~(\ref{umt}) we get then
\[
\p1(\E) \geq \p1(F) + \frac 12\, \sum_{i=1}^k P(E_i) - \xi\, k \p1(B)
\geq \bigg( V^{\frac{N-1}N}+ \frac {c_1 c_4}{8}\, k+\frac k 2-\xi k\bigg) \p1(B)\,.
\]
And finally, up to choosing $\xi<c_1c_4/8$, also in this case we get the validity of~(\ref{extraperimeter}), provided that $c_2\leq c_2''''' = c_1c_4/8-\xi$. Summarizing, we have a suitable choice of a small constant $\xi$, and in all the possible cases we have~(\ref{extraperimeter}), having chosen $c_2= \min\{c_2',\,c_2'',\, c_2''',\,c_2'''',\, c_2'''''\}$. The proof is then concluded.
\end{proof}

\section{Proof of Theorem~\mref{main}\label{secproof}}

This last section is devoted to the proof of Theorem~\mref{main}, which will be obtained by a careful use of the estimates proved in Section~\ref{sec:Euclest}. We start by considering a particular ambient space $\Omega$, that we can call ``double-universe''. The space is given by the disjoint union of two copies of $\R^N$, that we denote by $\R^N_1$ and $\R^N_2$. Both $\R^N_1$ and $\R^N_2$ are endowed with the standard Lebesgue volume, but with different notions of perimeter; more precisely, in $\R^N_1$ we consider the standard Euclidean perimeter (and we call it $P_1$ to help with the notation), while in $\R^N_2$ we consider the anisotropic perimeter corresponding to a smooth norm $\bar h$ (and we call it $P_2$). In other words, $\bar h:\R^N\to\R^+$ is an even, positively $1$-homogeneous function, smooth out of the origin and with strictly convex sublevels, and the perimeter of any set $F\subseteq \R^N_2$ is given by
\begin{equation}\label{defP2}
P_2(F) = \int_{\partial^* F} \bar h(\nu_F(x))\, d\H^{N-1}(x)\,.
\end{equation}
It is well-known that the unit volume isoperimetric set for this perimeter is the Wulff shape $W_{\bar h}$, which is the unit ball corresponding to the dual norm of $\bar h$, up to a homothety so to have unit volume. In particular, we pick a norm $\bar h$ such that the Wulff shape is almost exactly the unit cube $Q=[0,1]^N$ (it cannot be exactly the cube just because we want to build a smooth counterexample). Observe that, once we have decided what is the Wulff shape, we are still free to decide its perimeter; in fact, if we multiply the norm by a constant $\lambda>0$, the Wulff shape remains the same, but its perimeter is multiplied by $\lambda$. We can then fix $\bar h$ so that the perimeter $P_2$ of $W_{\bar h}$ is slightly larger than the perimeter $P_1$ of the ball $B$ with unit volume; more precisely, we assume that
\begin{equation}\label{QalmostB}
1 + \frac{c_2}2 < \frac{P_2(W_{\bar h})}{P_1(B)} < \frac{P_2(Q)}{P_1(B)} < 1+c_2\,,
\end{equation}
where $c_2$ is the constant in~(\ref{extraperimeter}) of Proposition~\ref{donow}.\par

Having defined volume and perimeter in $\R^N_1$ and in $\R^N_2$, we can do the same in the ``double universe'' $\Omega$. More precisely, since any subset $G$ of $\Omega$ is made by a part $E$ in $\R^N_1$ and a part $F$ in $\R^N_2$, we say that volume and perimeter of $G$ in $\Omega$ are simply the sum of the volumes and perimeters of $E$ and $F$. Analogously, any cluster $\G$ in $\Omega$ is made by a cluster $\E$ in $\R^N_1$ and a cluster $\F$ in $\R^N_2$, and volume and perimeter of $\G$ are defined as the sum of the volumes and perimeters of $\E$ and $\F$; we will call $P_\Omega$ the perimeter in $\Omega$. It is standard to observe that there exist minimal clusters of any volume in this ambient space. We are going to prove that for some volumes these minimal clusters have a non-empty part in each of the two subspaces $\R^N_1$ and $\R^N_2$.
\begin{prop}\label{propdoubleuniverse}
Let $\Omega$ be a double-universe as above, with a norm $\bar h$ such that~(\ref{QalmostB}) is in force, and let $V>0$ and $k\in\N$ be two sufficiently large constants, satisfying
\begin{align}\label{satisf1}
V^{\frac{N-1}N} > \frac{2k}{c_2} \,, && k> 2^{N(N+1)} \bigg(1+\frac 1{c_2}\bigg)^N\,.
\end{align}
Then, if $\G$ is a minimal $(k+1)$-cluster with volume $(V,\, 1,\, 1,\,\dots\,,\, 1)$, and $\E$ and $\F$ are its two parts in $\R^N_1$ and $\R^N_2$ respectively, both $\E$ and $\F$ are non-empty. The same holds if the volume of $\G$ is $t(V,\, 1,\, 1,\,\dots\,,\, 1)$, where $t>0$ is any constant.
\end{prop}
\begin{proof}
We can directly consider the case when the volume is $(V,\, 1,\, 1,\,\dots\,,\, 1)$, since the general case with the extra parameter $t$ is simply obtained by rescaling the clusters of a factor $t^{-1/N}$.\par

Let us assume that the claim is false, and let $\G$ be a cluster as in the claim, and such that either $\E$ or $\F$ is empty. We assume first that $\E=\emptyset$, so that the whole cluster is contained in $\R^N_2$. Then, calling $F_1$ the chamber of $\F$ with volume $V$, we have that
\[
P_\Omega(\G) = P_2(\F) \geq P_2(F_1) \geq V^{\frac{N-1}N} P_2(W_{\bar h})\,.
\]
We can take as competitor a cluster $\G'$ made by $k+1$ disjoint balls in $\R^N_1$, one having volume $V$ and the other $k$ having volume $1$. Then, we have
\[
P_\Omega(\G) \leq P_\Omega(\G') = P_1(B) \Big( V^{\frac{N-1}N} + k \Big)\,.
\]
Putting together the two above estimates and using~(\ref{QalmostB}), we get a contradiction with the left inequality in~\eqref{satisf1}.\par

Let us now assume that $\F=\emptyset$, so that the whole cluster is contained in $\R^N_1$. Then, by~(\ref{extraperimeter}) of Proposition~\ref{donow} we have
\begin{equation}\label{itsthere}
P_\Omega(\G) = P_1(\E) \geq \bigg( V^{\frac{N-1}N}+ \bigg( \frac 12 + c_2\bigg) k \bigg)\, P_1(B)\,.
\end{equation}
This time, we take as a competitor a cluster $\G''$ made by a Euclidean ball $E$ in $\R^N_1$ with volume $V$, and a cluster $\F''$ in $\R^N_2$ made by $k$ disjoint unit cubes, packed so to maximize the area of the common faces. If for instance $k=a^N$ for some $a\in\N$, then the union of the cubes is a cube of side $a$, and in that case the perimeter of $\F''$ is
\[
P_2(\F'') = \frac{a^N}2\, P_2(Q) + \frac{a^{N-1}}{2} P_2(Q) = \frac12\left(k\, + k^{\frac{N-1}N}\right)\, P_2(Q)\,.
\]
In general, if $k$ is not exactly an $N$-th power, a trivial estimate (which is not sharp at all, but enough for our purposes) gives
\[
P_2(\F'') \leq \bigg(\frac k 2 + 2^N k^{\frac{N-1}N}\bigg)\, P_2(Q)\,.
\]
Putting this estimate together with~(\ref{itsthere}), we get this time
\[\begin{split}
\bigg( V^{\frac{N-1}N}+ \bigg( \frac 12 + c_2\bigg) k \bigg)\, P_1(B) 
&\leq P_\Omega(\G) \leq P_\Omega(\G'') = P_1(E) + P_2(\F'')\\
&\leq V^{\frac{N-1}N} P_1(B)+ \bigg(\frac k 2\, + 2^N k^{\frac{N-1}N}\bigg)\, P_2(Q)\,,
\end{split}\]
which thanks to~(\ref{QalmostB}) immediately gives a contradiction with the right inequality in~(\ref{satisf1}).
\end{proof}

We are now able to give the proof of our main result.
\proofof{Theorem~\mref{main}}
We start defining $M$. The manifold $M$ is simply the $N$-dimensional sphere $\S^N\subseteq\R^{N+1}$ with the standard measure $\mu$. In order to properly define the Finsler structure $j_F$, we take three open sets $U_1\comp U_2\comp U_3$ contained in the north hemisphere $\{x_{N+1}>0\}$, and three open sets $U_4\comp U_5\comp U_6$ in the south one, we fix a norm $\bar h:\R^N\to\R^+$ with smooth and strictly convex Wulff shape $W_{\bar h}$ and satisfying~\eqref{QalmostB}, two constants $V$ and $\tilde k$ satisfying~\eqref{satisf1}, and we set $k=\tilde k+1$. Then, we select two charts $\varphi_3: \widetilde U_3 \to U_3$ and $\varphi_6:  \widetilde U_6 \to U_6$, where $\widetilde U_3$ and $\widetilde U_6$ are two open sets in $\R^N$, and we assume the charts to be \emph{isovolumetric}, that is, the push-forward ${\varphi_3^{-1}\!\!}_\# \,\mu$ (resp., ${\varphi_6^{-1}\!\!}_\#\, \mu$) coincides with the standard Lebesgue measure on $\widetilde U_3$ (resp., on $\widetilde U_6$).\par

We can then define a first smooth, reversible Finsler structure $\tilde j_F$ on $TM$ so that the following holds. If $\omega \in T_p^* M$ and $p\in U_3$, then we have $\tilde j_F^*(\omega) = |d(\varphi_3^{-1})_p^*(\omega)|$; instead, if $\omega\in T_p^* M$ and $p\in U_6$, then $\tilde j_F^*(\omega)=\bar h(d(\varphi_6^{-1})_p^*(-
\omega))$. For covectors $\omega\in T_p^* M$ with $p$ not belonging to $U_3$ or $U_6$ we have no requirement, except that the whole structure is smooth and reversible. Let us now call $\widetilde P_M$ the perimeter on $M$ induced by the structure $\tilde j_F$. By construction, keeping in mind formula~(\ref{weuse}) for the perimeter, whenever $E'\comp U_3$ and $E''\comp U_6$ we have
\begin{align}\label{Pervarphi36}
\widetilde P_M(E') = P\big(\varphi_3^{-1}(E')\big) \,, && \widetilde P_M(E'') = P_2\big(\varphi_6^{-1}(E'')\big)\,,
\end{align}
where as in Section~\ref{sec:Euclest} we denote by $P$ the standard Euclidean perimeter on $\R^N$, and as in~(\ref{defP2}) we denote by $P_2$ the perimeter induced by the norm $\bar h$.\par

The last ingredient that we need is a smooth function $\lambda: M \to (0,1]$ which satisfies $\lambda(p)=1$ if and only if $p\in\overline{U_1}\cup \overline{U_4}$; the final structure on $M$ will then be given by $j_F(v) = \lambda(p) \tilde j_F(v)$ for every $v\in T_p M$. We can now easily notice that lowering $\lambda$ has the effect of increasing the perimeter. In fact, if $j_F$ becomes smaller then so does the distance on $M$, and then the slopes of Lipschitz functions increase, and then so does the perimeter by the definition~(\ref{defP_M}). Alternatively, one can use~(\ref{weuse}) since lowering $j_F$ increases $j_F^*$. As a consequence, up to have chosen a function $\lambda$ which is sufficiently small outside of $U_2\cup U_5$, we can assume that for any set $E\subseteq M\setminus (U_2\cup U_5)$ one has
\begin{equation}\label{hugeP}
P_M(E) \geq 3 N\omega_N^{1/N} |E|_M^{\frac{N-1}N}\,,
\end{equation}
that is, the perimeter of $E$ is at least three times more than the Euclidean perimeter of a Euclidean ball having the same volume as $E$.\par

By construction, the triple $(M,j_F,\mu)$ is a compact, smooth, reversible Finsler manifold. Then, calling for brevity $w=\big( V,\, 1 ,\, 1 ,\, \dots \, , \, 1\big)\in(\R^+)^k=(\R^+)^{\tilde k+1}$, we define the sequence $\eps(j)$ as $\eps(j) = w / j$ for any $j$ large enough. Let then $\G$ be a minimal cluster in $M$ with volume $\eps(j)$, and assume by contradiction that $\G$ is connected. By~\cite[Theorem~A]{NP2025}, we know that the diameter of $\G$ is less than $K\cdot j^{-1/N}$ for some fixed constant $K$ depending on $k$, $V$, and on the structure of $M$.\par

We claim now that $\G$ is entirely contained either in $U_1$ or in $U_4$. Indeed, the diameter constraint implies that, as soon as $j$ is large enough, then $\G$ cannot contain points in $U_2\cup U_5$ together with points which are not in $U_3\cup U_6$. In other words, $\G$ is entirely contained either in $U_3$, or in $U_6$, or in $M\setminus (U_2\cup U_5)$. Suppose first that $\G\subseteq U_3$: then, up to a translation (via $\varphi_3$) we can move $\G$ to be entirely in $U_1$. Since $\lambda=1$ on $\overline{U_1}$ and $\lambda<1$ elsewhere, this translation strictly lowers the perimeter of the cluster unless $\G$ was already in $U_1$. In other words, if $\G$ is entirely contained in $U_3$ then it is actually entirely contained in $U_1$. In the very same way, if $\G$ is entirely contained in $U_6$ then it is actually entirely contained in $U_4$. Suppose now that $\G$ is contained in $M\setminus (U_2\cup U_5)$. Then, by~(\ref{hugeP}) we have
\[
P_M(\G) = \frac 12\, \sum_{i=0}^k P_M(G_i) 
\geq \frac 32\, \sum_{i=0}^k N\omega_N^{1/N} |G_i|_M^{\frac{N-1}N}
> \frac 32\, \sum_{i=1}^k N\omega_N^{1/N} |G_i|_M^{\frac{N-1}N}\,.
\]
However, we can take as competitor for $\G$ the cluster $\G'=\varphi_3(\widetilde \G)$, where $\widetilde\G=(\widetilde G_1,\, \widetilde G_2,\, \dots\, ,\, \widetilde G_k)$ is a cluster in $\varphi_3^{-1}(U_1)\subseteq \widetilde U_3\subseteq\R^N$ made by $k$ disjoint balls, one of which with volume $V/j$ and the other ones with volume $1/j$. Calling, for every $0\leq i\leq k$, $G'_i=\varphi_3(\widetilde G_i)$, keeping in mind~(\ref{Pervarphi36}), by minimality of $\G$ and by the above estimate, we have
\[\begin{split}
\frac 32\, \sum_{i=1}^k N\omega_N^{1/N} |G_i|_M^{\frac{N-1}N} &< P_M(\G)
\leq P_M(\G') = \frac 12\, \sum_{i=0}^k P_M(G_i') 
= \frac 12\, \sum_{i=0}^k \widetilde P_M(G_i') 
= \frac 12\, \sum_{i=0}^k P(\widetilde G_i) \\
&= \sum_{i=1}^k N\omega_N^{1/N} |\widetilde G_i|^{\frac{N-1}N}
= \sum_{i=1}^k N\omega_N^{1/N} |G_i|_M^{\frac{N-1}N}\,,
\end{split}\]
which is absurd. This concludes the proof of our claim, that is, $\G$ is entirely contained either in $U_1$ or in $U_4$.

Let us now consider the double-universe $\Omega$ of Proposition~\ref{propdoubleuniverse}: since the structure in $U_1$ is exactly the Euclidean one, and the structure in $U_4$ is exactly the one corresponding to the perimeter $P_2$ defined in~(\ref{defP2}), there is a function $\Phi:U_1 \cup U_4 \to \Omega$ which maintains the volume and the perimeter of any set or cluster. In particular, $\Phi(U_1)$ is an open set in $\R^N_1$, and $\Phi(U_4)$ is an open set in $\R^N_2$. The cluster $\Phi(\G)$ is then a cluster in $\Omega$ which is entirely contained either in $\R^N_1$ or in $\R^N_2$, and its volume is $\eps(j)=w/j$. Since the constants $V$ and $\tilde k$ satisfy~(\ref{satisf1}), Proposition~\ref{propdoubleuniverse} ensures that $\Phi(\G)$ is not a minimal cluster in $\Omega$. Thus, there is a cluster $\widehat \G$ in $\Omega$ with perimeter strictly smaller than that of $\G$. Thanks to~\cite[Theorem~A]{NP2025}, we know that $\widehat \G$ is made by at most $k$ connected components, each one with diameter of order $j^{-1/N}$. Thus, as soon as $j$ is large enough, up to a translation $\widehat\G$ is contained in $\Phi(U_1\cup U_4)$. But then, the cluster $\Phi^{-1}(\widehat \G)$ is a cluster in $M$ with volume $\eps(j)$ and strictly less perimeter than $\G$, which is against the assumptions. Summarizing, we have proved the non-connectedness of small minimal clusters in $M$.\par

To conclude the proof, we have to show the non-existence of a norm $\|\cdot\|$ on $\R^N$ such that~(\ref{same-limit}) holds. Thus, we assume that $\|\cdot\|$ is such a norm and we look for a contradiction. The key idea is that the minimal clusters in $M$ with mass $\eps(j)$ are all equal, up to rescaling. Let us be more precise. Let $j$ be a sufficiently large integer, and let $\G$ be a minimal cluster in $M$ with volume $\eps(j)$. We already know that $\G$ is not connected and that we can write $\G$ as the union of two clusters $\E$ and $\F$, contained in $U_1$ and $U_4$ respectively. We define then $w'= j |\E|$ and $w''= j|\F|$, so that $w=w'+w''$. Again by the small diameter property, we deduce that both $\E$ and $\F$ are minimal clusters in $M$, and thus
\[
\J_{M,k} (w/j) = P_M(\G) = P_M(\E) + P_M(\F) = \J_{M,k}(w'/j) + \J_{M,k}(w''/j)\,.
\]
Now, $\Phi(\G)$ is a cluster in $\Omega$ with volume $\eps(j)$. Since we can perform rescaling in $\Omega$, for any $j'$ large enough we can consider
\[
\widehat \G' = \bigg(\frac j{j'}\bigg)^{\frac 1N}\, \Big( \Phi(\G)\Big)\,,
\]
which is a cluster in $\Omega$ with volume $\eps(j')$. As noticed above, as soon as $j$ and $j'$ are large enough we have that $\widehat\G'$ is contained in $\Phi(U_1\cup U_4)$, so we can set $\G'=\Phi^{-1}(\widehat\G')$, which is a cluster in $M$ with volume $\eps(j')$. Since the map $\Phi$ keeps both volumes and perimeters unchanged, we have
\[
P_M(\G' ) = \bigg(\frac j{j'}\bigg)^{\frac {N-1}N}\, P_M(\G) \,,
\]
which implies
\[
\J_{M,k} ( w/j') \leq  \bigg(\frac j{j'}\bigg)^{\frac {N-1}N}\, \J_{M,k}(w/j)\,.
\]
Exchanging the roles of $j$ and $j'$, we get that the above inequality is in fact an equality. This means that for every $j'$ large enough we have
\[
\J_{M,k} ( w/j') j'^{\frac {N-1}N} = P_M(\G)\, j^{\frac {N-1}N}\,,
\]
and similarly
\begin{align*}
\J_{M,k} ( w'/j') j'^{\frac {N-1}N} = P_M(\E)\, j^{\frac {N-1}N}\,, && \J_{M,k} ( w''/j') j'^{\frac {N-1}N} = P_M(\F)\, j^{\frac {N-1}N}\,.
\end{align*}
By the rescaling properties of the perimeter on $\R^N$ with a fixed norm we know that
\[
\J_{\R^N,\|\cdot\|,k}(w) = \J_{\R^N,\|\cdot\|,k}(w/j') j'^{\frac{N-1}N}
\]
for every $j'\in\N$, and the same equality holds with $w'$ or $w''$ in place of $w$. But then, since we are assuming the validity of~(\ref{same-limit}), we obtain 
\[\begin{split}
\J_{\R^N,\|\cdot\|,k}(w) &= \lim_{j'\to +\infty} \J_{\R^N,\|\cdot\|,k}(w/j') j'^{\frac{N-1}N}
= \lim_{j'\to +\infty} \frac{\J_{\R^N,\|\cdot\|,k}(w/j') }{\J_{M,k}(w/j')}\, \J_{M,k}(w/j') j'^{\frac{N-1}N}\\
&= \lim_{j'\to +\infty}  \J_{M,k}(w/j') j'^{\frac{N-1}N}=P_M(\G)\, j^{\frac {N-1}N}\,,
\end{split}\]
and in the very same way we have also
\begin{align*}
\J_{\R^N,\|\cdot\|,k}(w') = P_M(\E)\, j^{\frac {N-1}N}\,, && \J_{\R^N,\|\cdot\|,k}(w'') = P_M(\F)\, j^{\frac {N-1}N}\,.
\end{align*}
Since $P_M(\G)=P_M(\E)+P_M(\F)$, we deduce
\[
\J_{\R^N,\|\cdot\|,k}(w) = \J_{\R^N,\|\cdot\|,k}(w')+\J_{\R^N,\|\cdot\|,k}(w'')\,.
\]
And finally, since both $w'$ and $w''$ are not null, it is well-known that this equality is impossible, because minimal clusters in $\R^N$ with a fixed norm are always connected. Indeed, if the equality was true, then the disjoint union of two minimal clusters with volumes $w'$ and $w''$ respectively would be a non-connected minimal cluster with volume $w$.
\end{proof}

\bigskip
\subsection*{Acknowledgments}
This research was funded in part by the Austrian Science Fund (FWF) [grant DOI \href{https://www.fwf.ac.at/en/research-radar/10.55776/EFP6}{10.55776/EFP6}]. For open access purposes, the authors have applied a CC BY public copyright license to any author-accepted manuscript version arising from this submission. Both authors are member of the Istituto Nazionale di Alta Matematica (INdAM), Gruppo Nazionale per l'Analisi Matematica, la Probabilit\`a e le loro Applicazioni (GNAMPA). D.C. was partially supported by the INdAM--GNAMPA 2024 Project \textit{Propriet\`a geometriche e regolarit\`a in problemi variazionali locali e non locali}, code CUP\_E53\-C23\-001\-670\-001; A.P. was partially supported by the PRIN Project 2022E9CF89, ``Geometric Evolution Problems and Shape Optimization (GEPSO)'', financed by the European Union, and by the DFG through the Grant IRTG3132 ``Analysis of Nonlinear PDE'' (Germany).

\bibliographystyle{plain}

\begin{thebibliography}{999}
\bibitem{A2025} G. Antonelli, Isoperimetric problems and lower bounds on curvature, Anal. Geom. Metr. Spaces {\bf 14} (2026), no. 1.                    
\bibitem{ABFP2022} G. Antonelli, E. Bru\`e, M. Fogagnolo \& M. Pozzetta, On the existence of isoperimetric regions in manifolds with nonnegative Ricci curvature and Euclidean volume growth, Calc. Var. Partial Differential Equations {\bf 61} (2022), no. 2, Paper No. 77.
\bibitem{AG2024} G. Antonelli \& F. Glaudo, Nonexistence of isoperimetric sets in spaces of positive curvature, J. Reine Angew. Math. {\bf 812} (2024), 211--228.
\bibitem{APPS2025} G. Antonelli, E. Pasqualetto, M. Pozzetta \& D. Semola, Sharp isoperimetric comparison on non-collapsed spaces with lower Ricci bounds, Ann. Sci. \'Ec. Norm. Sup\'er. (4) {\bf 58} (2025), no. 1, 1--52.
\bibitem{BCS2000} D. Bao, S.-S. Chern \& Z. Shen, An introduction to Riemann-Finsler geometry, Graduate Texts in Mathematics, volume 200, Springer-Verlag, 2000.
\bibitem{CM2022} F. Cavalletti \& D. Manini, Isoperimetric inequality in noncompact MCP spaces, Proc. Amer. Math. Soc. {\bf 150} (2022), no. 8, 3537--3548.
\bibitem{CM2024} F. Cavalletti \& D. Manini, Rigidities of Isoperimetric inequality under nonnegative Ricci curvature, J. Eur. Math. Soc. {\bf 28} (2026), no. 8, 3235--3294.
\bibitem{DeRN2025} A. De Rosa \& R. Neumayer, Local minimizers of the anisotropic isoperimetric problem on closed manifolds, Indiana Univ. Math. J. {\bf 74} (2025), no. 5, 1423--1466.
\bibitem{FPP} L. Felicetti, R. Ponciano \& A. Pratelli, Connectedness of small minimal clusters in fixed-norm Finsler manifolds, preprint (2026).
\bibitem{FM2011} A. Figalli \& F. Maggi, On the shape of liquid drops and crystals in the small mass regime, Arch. Ration. Mech. Anal. {\bf 201} (2011), no. 1, 143--207.
\bibitem{FJ} N. Fusco \& V. Julin, A strong form of the quantitative isoperimetric inequality, Calc. Var. Partial Differential Equations {\bf 50} (2014), no. 3--4, 925--937.
\bibitem{manini2024} D. Manini, Isoperimetric inequality for Finsler manifolds with non-negative Ricci curvature, Rev. Mat. Iberoam. {\bf 40} (2024), no. 5, 1631--1690.
\bibitem{MJ2000} F. Morgan \& D.L.\, Johnson, Some sharp isoperimetric theorems for Riemannian manifolds, Indiana Univ. Math. J. {\bf 49} (2000), no. 3, 1017--1041.
\bibitem{N2009} S. Nardulli, The isoperimetric profile of a smooth Riemannian manifold for small volumes, Ann. Global Anal. Geom. {\bf 36} (2009), no. 2, 111--131.
\bibitem{NP2025} S. Nardulli \& A. Pratelli, Connectedness properties of small minimal clusters in Riemannian or Finsler manifolds, preprint (2025).
\bibitem{ohta2009} S. Ohta, Finsler interpolation inequalities, Calc. Var. Partial Differential Equations {\bf 36} (2009), no. 2, 211--249.
\bibitem{ohta2018} S. Ohta, Needle decompositions and isoperimetric inequalities in Finsler geometry, J. Math. Soc. Japan {\bf 70} (2018), no. 2, 651--693.
\bibitem{Ohtabook} S. Ohta, Comparison Finsler geometry, Springer Monographs in Mathematics, 2021.
\bibitem{P2023} M. Pozzetta, Isoperimetry on manifolds with Ricci bounded below: overview of recent results and methods, in: Anisotropic Isoperimetric Problems and Related Topics, Springer INdAM Ser (2024), 49--89.
\end{thebibliography}

\end{document}